# Random sets of isomorphism of linear operators on Hilbert space

Roman Vershynin[1],*

*University of California, Davis*

**Abstract:** This note deals with a problem of the probabilistic Ramsey theory in functional analysis. Given a linear operator $T$ on a Hilbert space with an orthogonal basis, we define the *isomorphic structure* $\Sigma(T)$ as the family of all subsets of the basis so that $T$ restricted to their span is a nice isomorphism. Our main result is a dimension-free optimal estimate of the size of $\Sigma(T)$. It improves and extends in several ways the principle of restricted invertibility due to Bourgain and Tzafriri. With an appropriate notion of randomness, we obtain a randomized principle of restricted invertibility.

## 1. Introduction

### 1.1. Randomized Ramsey-type problems

Finding a nice structure in a big unstructured object is a recurrent theme in mathematics. This direction of thought is often called Ramsey theory, although Ramsey theory was originally only associated with combinatorics. One celebrated example is Van der Waerden's theorem: for any partition of the integers into two sets, one of these sets contains arbitrary long arithmetic progressions.

Ramsey theory meets probability theory when one asks about the quality of *most* sub-structures of a given structure. Can one improve the quality of a structure by passing to its *random* sub-structure? (a random subgraph, for example). A remarkable example of the randomized Ramsey theory is Dvoretzky's theorem in geometric functional analysis in the form of V. Milman (see [4], 4.2). One of its corollaries states that, for any $n$-dimensional finite-dimensional Banach space, a random $O(\log n)$-dimensional subspace (with respect to some natural measure) is well isomorphic to a Hilbert space.

### 1.2. The isomorphism structure of a linear operator

In this note we are trying to find a nice structure in an arbitrary bounded linear operator on a separable Hilbert space. Let $T$ be a bounded linear operator on a Hilbert space $H$ with an orthonormal basis $(e_i)_{i\in\mathbb{N}}$. We naturally think of $T$ as being nice if it is a nice isomorphism on $H$. However, this situation is rather rare; instead, $T$ may be a nice isomorphism on the subspace spanned by some subsets of the basis. So, instead of being a "global" isomorphism, $T$ may be a "local" isomorphism when restricted to certain subspaces of $H$. A central question is then – how many such subspaces are there? Let us call these subspaces an isomorphism structure of $T$:

---

[1]Department of Mathematics, University of California, Davis, CA 95616, USA.
*The author is an Alfred P. Sloan Research Fellow. He was also supported by NSF DMS 0401032.

*AMS 2000 subject classifications:* 46B09, 47D25.

*Keywords and phrases:* restricted invertibility, operators on Hilbert spaces.

**Definition 1.1.** Let $T$ be a bounded linear operator on a Hilbert space $H$, and $(e_i)_{i \in \mathbb{N}}$ be an orthonormal basis of $H$. Let $0 < \varepsilon < 1$. A set $\sigma$ of $\mathbb{N}$ is called a *set of $\varepsilon$-isomorphism* of $T$ if the equivalence

$$(1) \qquad (1-\varepsilon) \sum_{i \in \sigma} \|a_i T e_i\|^2 \leq \|\sum_{i \in \sigma} a_i T e_i\|^2 \leq (1+\varepsilon) \sum_{i \in \sigma} \|a_i T e_i\|^2$$

holds for every choice of scalars $(a_i)_{i \in \sigma}$. The *$\varepsilon$-isomorphism structure* $\Sigma(T, \varepsilon)$ consists of all such sets $\sigma$.

How big is the isomorphism structure? From the probabilistic point of view, we can ask for the probability that a random subset of (a finite interval of) the basis is the set of isomorphism. Unfortunately, this probability is in general exponentially small. For example, if $T$ acts as $Te_i = e_{\lceil (i+1)/2 \rceil}$, then every set of isomorphism contains no pairs of the form $\{2i-1, 2i\}$. Hence a random subset of a finite interval is unlikely to be a set of isomorphism of $T$. However, an appropriate notion of randomness yields a clean optimal bound on the size of the isomorphic structure. This is the main result of this note, which extends in several ways the Bourgain-Tzafriri's principle of the restricted invertibility [1], as we will see shortly.

**Theorem 1.2.** *Let $T$ be a norm-one linear operator on a Hilbert space $H$, and let $0 < \varepsilon < 1$. Then there exists a probability measure $\nu$ on the isomorphism structure $\Sigma(T, \varepsilon)$, such that*

$$(2) \qquad \nu\{\sigma \in \Sigma(T, \varepsilon) \mid i \in \sigma\} \geq c\varepsilon^2 \|Te_i\|^2 \quad \text{for all } i.$$

Here and thereafter $c, C, c_1, \ldots$ denote positive absolute constants.

Theorem 1.2 gives a lower bound on the average of the characteristic functions of the sets of the isomorphism. Indeed, the left hand side in (2) clearly equals $\int_{\Sigma(T,\varepsilon)} \chi_\sigma(i) \, d\nu(\sigma)$. Thus, in absence of "true" randomness in the isomorphic structure $\Sigma(T, \varepsilon)$, we can still measure the size of $\Sigma(T, \varepsilon)$ by bounding below the *average of the characteristic functions* of its sets. It might be that considering this weak type of randomness might help in other problems, in which the usual, strong randomness, fails.

### 1.3. Principle of restricted invertibility

One important consequence of Theorem 1.2 is that *there always exists a big set of isomorphism of $T$*. This extends and strengthens a well known result due to Bourgain and Tzafriri, known under the name of the principle of restricted invertibility [1]. We will show how to find a big set of isomorphism; its size can be measured with respect to an arbitrary measure $\mu$ on $\mathbb{N}$. For the rest of the paper, we denote the measure of the singletons $\mu(\{i\})$ by $\mu_i$. Summing over $i$ with weights $\mu_i$ in (2) and using Theorem 1.2, we obtain

$$(3) \qquad \begin{aligned} \int_{\Sigma(T,\varepsilon)} \mu(\sigma) \, d\nu(\sigma) &= \sum_i \mu_i \int_{\Sigma(T,\varepsilon)} \chi_\sigma(i) \, d\nu(\sigma) \\ &= \sum_i \mu_i \, \nu\{\sigma \in \Sigma(T,\varepsilon) \mid i \in \sigma\} \geq c\varepsilon^2 \sum_i \mu(i) \|Te_i\|^2. \end{aligned}$$

Replacing the integral in the left hand side of (4) by the maximum shows that there exists a big set of isomorphism:

**Corollary 1.3.** *Let $T$ be a norm-one linear operator on a Hilbert space $H$, and let $\mu$ be a measure on $\mathbb{N}$. Then, for every $0 < \varepsilon < 1$, there exists a set of $\varepsilon$-isomorphism $\sigma$ of $T$ such that*

$$\mu(\sigma) \geq c\varepsilon^2 \sum_i \mu_i \|Te_i\|^2. \tag{4}$$

Earlier, Bourgain and Tzafriri [1] proved a weaker form of Corollary 1.3 with only the *lower* bound in the definition (1) of the set of isomorphism, for a uniform measure $\mu$ on an interval, under an additional assumption on the uniform lower bound on $\|Te_i\|$, and for some fixed $\varepsilon$.

**Theorem 1.4** (Bourgain-Tzafriri's principle of restricted invertibility). *Let $T$ be a linear operator on an $n$-dimensional Hilbert space $H$ with an orthonormal basis $(e_i)$. Assume that $\|Te_i\| = 1$ for all $i$. Then there exits a subset $\sigma$ of $\{1, \ldots, n\}$ such that $|\sigma| \geq cn/\|T\|^2$ and*

$$\|Tf\| \geq c\|f\|$$

*for all $f \in \mathrm{span}(e_i)_{i \in \sigma}$.*

This important result has found applications in Banach space theory and harmonic analysis. Corollary 1.3 immediately yields a stronger result, which is dimension-free and which yields an almost isometry:

**Corollary 1.5.** *Let $T$ be a linear operator on a Hilbert space $H$ with an orthonormal basis $(e_i)$. Assume that $\|Te_i\| = 1$ for all $i$. Let $\mu$ be a measure on $\mathbb{N}$. Then, for every $0 < \varepsilon < 1$, there exits a subset $\sigma$ of $\mathbb{N}$ such that $\mu(\sigma) \geq c\varepsilon^2/\|T\|^2$ and such that*

$$(1-\varepsilon)\|f\| \leq \|Tf\| \leq (1+\varepsilon)\|f\| \tag{5}$$

*for all $f \in \mathrm{span}(e_i)_{i \in \sigma}$.*

Szarek [5] proved a weaker form of Corollary 1.3 with only the *upper* bound in the definition (1) of the set of isomorphism, and with some fixed $\varepsilon$.

For the counting measure on $\mathbb{N}$, Corollary 1.3 was proved in [7]. In this case, bound (4) reads as

$$|\sigma| \geq c\varepsilon^2 \|T\|_{\mathrm{HS}}^2, \tag{6}$$

where $\|T\|_{\mathrm{HS}}$ denotes the Hilbert-Schmidt norm of $T$. (If $T$ is not a Hilbert-Schmidt operator, then an infinite $\sigma$ exists).

## 2. Proof of Theorem 1.2

Corollary 1.3 is a consequence of two suppression results due to Szarek [5] and Bourgain-Tzafriri [2]. We will then deduce Theorem 1.2 from Corollary 1.3 by a simple separation argument from [2].

To prove Corollary 1.3, we can assume by a straighforward approximation that our Hilbert space $H$ is finite dimensional. We can thus identify $H$ with the $n$-dimensional Euclidean space $\ell_2^n$, and identify the basis $(e_i)_{i=1}^n$ of $H$ with the canonical basis of $\ell_2^n$. Given a subset $\sigma$ of $\{1, \ldots, n\}$ (or of $\mathbb{N}$), by $\ell_2^\sigma$ we denote the subspace of $\ell_2^n$ (of $\ell_2$ respectively) spanned by $(e_i)_{i \in \sigma}$. The orthogonal projection onto $\ell_2^\sigma$ is denoted by $Q_\sigma$.

With a motivaiton different from ours, Szarek proved in ([5], Lemma 4) the following suppression result for operators in $\ell_2^n$.

**Theorem 2.1** (Szarek). *Let $T$ be a norm-one linear operator on $\ell_2^n$. Let $\lambda_1, \ldots, \lambda_n$, $\sum_{i=1}^n \lambda_i = 1$, be positive weights. Then there exists a subset $\sigma$ of $\{1, \ldots, n\}$ such that*

$$\sum_{i \in \sigma} \lambda_i \|Te_i\|^{-2} \geq c \tag{7}$$

*and such that the inequality*

$$\|\sum_{i \in \sigma} a_i Te_i\|^2 \leq C \sum_{i \in \sigma} \|a_i Te_i\|^2$$

*holds for every choice of scalars $(a_i)_{i \in \sigma}$.*

**Remark 2.2.** Inequality (7) for a probability measure $\lambda$ on $\{1, \ldots, n\}$ is equivalent to the inequality

$$\mu(\sigma) \geq c \sum_i \mu_i \|Te_i\|^2 \tag{8}$$

for a positive measure $\mu$ on $\{1, \ldots, n\}$.

Indeed, (7) implies (8) with

$$\lambda_i = \frac{\mu_i \|Te_i\|^2}{\sum_i \mu_i \|Te_i\|^2}.$$

Conversely, (8) implies (7) with $\mu_i = \lambda_i \|Te_i\|^{-2}$.

Theorem 2.1 and Remark 2.2 yield a weaker version of Corollary 1.3 – with only the upper bound in the definition (1) of the set of isomorphism, and with some fixed $\varepsilon$.

To prove Corollary 1.3 in full strength, we will use the following suppression analog of Theorem 1.2 due to Bourgain and Tzafriri [2].

**Theorem 2.3** (Bourgain-Tzafriri). *Let $S$ be a linear operator on $\ell_2$ whose matrix relative to the unit vector basis has zero diagonal. For a $\delta > 0$, denote by $\Sigma'(S, \delta)$ the family of all subsets $\sigma$ of $\mathbb{N}$ such that $\|Q_\sigma S Q_\sigma\| \leq \delta \|S\|$. Then there exists a probability measure $\nu'$ on $\Sigma'(S, \delta)$ such that*

$$\nu'\{\sigma \in \Sigma'(S, \delta) \mid i \in \sigma\} \geq c\delta^2 \quad \text{for all } i. \tag{9}$$

*Proof of Corollary 1.3.* We define a linear operator $T_1$ on $H = \ell_2^n$ as

$$T_1 e_i = Te_i / \|Te_i\|, \quad i = 1, \ldots, n.$$

Theorem 2.1 and the remark below it yield the existence of a subset $\sigma$ of $\{1, \ldots, n\}$ whose measure satisfies (8) and such that the inequality

$$\|T_1 f\| \leq C \|f\|$$

holds for all $f \in \text{span}(e_i)_{i \in \sigma}$. In other words, the operator

$$T_2 = T_1 Q_\sigma$$

satisfies

(10) $$\|T_2\| \leq C.$$

We will apply Theorem 2.3 for the operator $S$ on $\ell_2^\sigma$ defined as

(11) $$S = T_2^* T_2 - I \quad \text{and with} \quad \delta = \varepsilon/\|S\|.$$

Indeed, $S$ has zero diagonal:

$$\langle Se_i, e_i \rangle = \|T_2 e_i\|^2 - 1 = \|T_1 e_i\|^2 - 1 = 0 \quad \text{for all } i \in \sigma.$$

Also, $S$ has nicely bounded norm by (10):

$$\|S\| \leq \|T_2\|^2 + 1 \leq C^2 + 1,$$

which yields a lower bound on $\delta$:

(12) $$\delta \geq \varepsilon/(C^2 + 1).$$

So, Theorem 2.3 yields a family $\Sigma'(S,\delta)$ of subsets of $\sigma$ and a measure $\nu'$ on this family. It follows as before that $\Sigma'(S,\delta)$ must contain a big set, because

$$\int_{\Sigma'(S,\delta)} \mu(\sigma') \, d\nu'(\sigma') = \sum_{i \in \sigma} \mu_i \int_{\Sigma'(S,\delta)} \chi_{\sigma'}(i) \, d\nu'(\sigma')$$
$$= \sum_{i \in \sigma} \mu_i \, \nu'\{\sigma' \in \Sigma'(S,\delta) \mid i \in \sigma'\}$$
$$\geq \sum_{i \in \sigma} \mu_i \cdot c\delta^2 \geq c'\varepsilon^2 \mu(\sigma)$$

where the last inequality follows from (12) with $c' = c(C^2 + 1)^{-2}$. Thus there exists a set $\sigma' \in \Sigma'(S,\delta)$ such that by (8) we have

$$\mu(\sigma') \geq c'\varepsilon^2 \mu(\sigma) \geq c''\varepsilon^2 \sum_{i=1}^n \mu_i \|Te_i\|^2,$$

so with the measure as required in (4).

It remains to check that $\sigma'$ is a set of $\varepsilon$-isomorphism of $T$. Consider an $f \in \text{span}(e_i)_{i \in \sigma'}$, $\|f\| = 1$. By the suppression estimate in Theorem 2.3 and by our choice of $S$ and $\delta$ made in (11), we have

$$\varepsilon = \delta\|S\| \geq |\langle Q_{\sigma'} S Q_{\sigma'} f, f \rangle|$$
$$= |\langle Sf, f \rangle| \quad \text{because } Q_{\sigma'} f = f$$
$$= |\|T_2 f\|^2 - \|f\|^2| \quad \text{by the definition of } S$$
$$= |\|T_1 f\|^2 - 1| \quad \text{because } Q_{\sigma'} f = Q_\sigma f = f \text{ as } \sigma' \subset \sigma.$$

It follows by homogeneity that

$$(1-\varepsilon)\|f\|^2 \leq \|T_1 f\|^2 \leq (1+\varepsilon)\|f\|^2 \quad \text{for all } f \in \text{span}(e_i)_{i \in \sigma'}.$$

By the definition of $T_1$, this means that $\sigma'$ is a set of $\varepsilon$-isomorphism of $T$. This completes the proof. $\square$

*Proof of Theorem 1.2.* We deduce Theorem 1.2 from Corollary 1.3 by a separation argument, which is a minor adaptation of the proof of Corollary 1.4 in [2].

We first note that, by Remark 2.2, an equivalent form of the consequence of Corollary 1.3 is the following. For every probability measure $\lambda$ on $\mathbb{N}$, there exists a set $\sigma \in \Sigma(T,\varepsilon)$ such that

$$\sum_{i \in \sigma} \lambda_i \|Te_i\|^{-2} \geq c\varepsilon^2. \tag{13}$$

We consider the space of continuous functions $C(\Sigma(T,\varepsilon))$ on the isomorphism structure $\Sigma(T,\varepsilon)$, which is compact in its natural topology (of pointwise convergence of the indicators of the sets $\sigma \in \Sigma(T,\varepsilon)$). For each $i \in \mathbb{N}$, define a function $\pi_i \in C(\Sigma(T,\varepsilon))$ by setting

$$\pi_i(\sigma) = \chi_\sigma(i)\, \|Te_i\|^{-2}, \quad \sigma \in \Sigma(T,\varepsilon).$$

Let $\mathcal{C}$ be the convex hull of the set of functions $\{\pi_i, i \in \mathbb{N}\}$. Every $\pi \in \mathcal{C}$ can be expressed a convex combination $\pi = \sum_i \lambda_i \pi_i$. By Corollary 1.3 in the form (13), there exists a set $\sigma \in \Sigma(T,\varepsilon)$ such that $\pi(\sigma) \geq c\varepsilon^2$. Thus $\|\pi\|_{C(\Sigma(T,\varepsilon))} \geq c\varepsilon^2$. We conclude by the Hahn-Banach theorem that there exists a probability measure $\nu \in C(\Sigma(T,\varepsilon))^*$ such that

$$\nu(\pi) = \int_{\Sigma(T,\varepsilon)} \pi(\sigma)\, d\nu(\sigma) \geq c\varepsilon^2 \quad \text{for all } \pi \in \mathcal{C}.$$

Applying this estimate for $\pi = \pi_i$, we obtain

$$\int_{\Sigma(T,\varepsilon)} \chi_\sigma(i)\, d\nu(\sigma) \geq c\varepsilon^2 \|Te_i\|^2,$$

which is exactly the conclusion of the theorem. □

**Remark 2.4.** The proof of Theorem 1.2 given above is a combination of previously known tools – two suppression results due to [5] and [2] and a separation argument from [2]. The new point was to realize that the suppression result of Szarek [5], developed with a different purpose in mind, gives a sharp estimate when combined with the results of [2]. To find a set of the isomorphism as in (1), one needs to reduce the norm of the operator with [5] *before* applying restricted invertibility principles from [2].

**Acknowledgement.** I am grateful to the referee for numerous comments and suggestions.